\documentclass[12pt]{article}
\usepackage{amsmath, amssymb}
\usepackage{amsfonts}
\usepackage{theorem}
\usepackage{makeidx}

{\theorembodyfont{\itshape}}

\newtheorem{theorem}{Theorem} 

\newtheorem{lemma}{Lemma}[section]

\numberwithin{equation}{section}

\newcommand{\T}{\, \mathbf{T} \, }
\newcommand{\Res}{\, \mathbf{R} \, }

\newcommand{\R}{\mathbb{R}}
\newcommand{\Z}{\mathbb{Z}}

\newcommand{\e}{\varepsilon}

\newcommand{\myproof}{\noindent {\bf Proof}: \quad}

\newcommand{\myendproof}{\hspace*{\fill}{{$\square$}} \vspace{10pt}}

\newcommand{\donothing}[1]{{}}

\newcommand{\Av}{\operatorname*{Av}}
\newcommand{\xb}{\mathbf{x}}
\newcommand{\yb}{\mathbf{y}}
\newcommand{\pb}{\mathbf{p}}

\newcommand{\rb}{\mathbf{r}}
\renewcommand{\sb}{\mathbf{s}}

\newcommand{\E}{{\cal E}}

\newcommand{\A}{{\mathcal A}}
\newcommand{\calS}{{\mathcal S}}
\newcommand{\calL}{{\mathcal L}}

\newcommand{\no}{\nonumber}

\donothing{
\newcommand{\xb}{{x}}
\newcommand{\yb}{{y}}
\newcommand{\pb}{{p}}

}

\renewcommand{\l}{\langle \!\langle}
\renewcommand{\r}{\rangle \!\rangle}
\newcommand{\la}{\langle}
\newcommand{\ra}{\rangle}

\begin{document}

\title{Superdiffusivity  of  asymmetric
 exclusion process in dimensions one and two}

\author{C. Landim\footnote{Partially supported by CNRS, UMR 6085, landim@impa.br},$\;$
J. Quastel\footnote{Partially supported by
NSERC, quastel@math.toronto.edu},$\;$
M. Salmhofer\footnote{mns@mis.mpg.de},$\;$ and
H.-T. Yau\footnote{Partially supported by NSF grant
DMS-0072098,  yau@cims.nyu.edu} \\[1ex]{\footnotesize
 IMPA,  University of Toronto, 
Max--Planck--Institut, Leipzig, and Courant Institute} }




\maketitle

\begin{abstract}
We prove that the diffusion coefficient for
the asymmetric exclusion process
diverges at least as fast as $t^{1/4}$ in dimension $d=1$ and 
$(\log t)^{1/2}$ in $d=2$. The method
applies to nearest and non-nearest neighbor
asymmetric exclusion processes.
\end{abstract}


\section{Introduction}

Asymmetric exclusion is a
Markov process on $\{0,1\}^{\mathbb Z^d}$ consisting of interacting
continuous time random walks with asymmetric jump rates.
There is at most one particle allowed per site.
A particle at a site $x$ waits for an exponential time and then
jumps to $y$ provided the site is not occupied.
Otherwise the jump is suppressed  and the process starts again.  The jump
is attempted at rate $p(y-x)$.
In this article the jump law $p(\cdot)$ is assumed to have a nonzero mean,
 so that there is 
transport of the system. 

Consider the system in equilibrium with a Bernoulli product measure of
density $\rho$ as the invariant measure. Define the time dependent correlation 
function in equilibrium by
$$
S(x, t)=
\langle  \eta_x(t);
\eta_0(0) \rangle
$$
If we choose $\rho=1/2$,
there is no net global drift, i.e., $\sum_x x S(x, t) = 0$.
Otherwise one needs to  subtract a net drift, which complicates but does 
not change the results or methods.   Our main question is the behaviour
for large $t$ of the diffusion coefficient,
$$
D(t) = \frac{1}{ 4t} \sum_x x^2 S(x, t).
$$

In dimensions $d \ge 3$, the diffusion coefficient was proved to be bounded
 \cite{LY} for general asymmetric simple exclusion processes.
Based on mode coupling theory,
Beijeren, Kutner and  Spohn \cite{Bei Ku Spo}
conjectured that $D(t)\sim
(\log t)^{2/3}$ in dimension $d=2$ and
$D(t)\sim t^{1/3}$ in $d=1$.  Similar predictions we made in 
\cite{KPZ} for the Kardar-Parisi-Zhang equation in $d=1$, which, when
differentiated and appropriately discretized, yields the asymmetric 
exclusion process.

This problem has received much attention recently in the context
of integrable systems. The main quantity analyzed there
is fluctuation of the current
across the origin in
$d=1$ for  the
totally asymmetric simple exclusion process (only nearest neighbor jumps to 
the right), starting from the special initial configuration 
with all sites to the 
left of the origin
occupied and all sites to the right
of the origin empty.
Johansson \cite{Joh} observed that in this special situation
the current across the origin
can be mapped
into a last passage percolation problem.
By analyzing  this 
problem asymptotically, Johansson proved that the variance of
the current is of order  $t^{2/3}$.
In the case of discrete time,
Baik and Rains \cite{BR} analyze an extended version of the last passage
percolation problem and  obtain fluctuations of order $t^\alpha$,
where $\alpha =1/3$ or $\alpha =1/2$ depending on the parameters of
the model.
Both the approaches of  \cite{Joh} and
\cite{BR} are related to the earlier results of
Baik-Deift-Johansson \cite{Ba dei jo} on the distribution
of the length of the longest increasing subsequence in random
permutations.

In \cite{PS} (see also \cite{Prae Spo2}),  Pr\"{a}hofer and Spohn
succeeded  in mapping
the current of the totally asymmetric simple exclusion process
 into a last passage percolation problem
for a general class of initial
data, including the equilibrium case considered in this article.
For the discrete time case, the extended problem is
closely related to the work
\cite{BR}, but the boundary conditions
are  different. For continuous time, besides the
boundary condition issue, one needs to  extend
the result of \cite{BR}
from the geometric to
the exponential distribution.

To relate these results to our problem,
the variance of the current across the origin
is proportional to
\begin{equation}\label{0.1}
\sum_x |x| S(x, t)
\end{equation}
Therefore, Johansson's result on the variance of
the current  can be interpreted
as the spreading of  $S(x,t)$ being of order $t^{2/3}$.
If we combine the work of \cite{PS}
and \cite{BR}, neglect various issues discussed above, and extrapolate
to the second moment, we obtain a growth of the second moment as
$t^{4/3}$, consistent with the conjectured $D(t) \sim t^{1/3}$.

The results based on integrable
systems give not just the variance of the current
across the origin, but also its limiting distribution.
The main restrictions
appear to be the rigid requirements on the fine details of the dynamics
and on the initial data, the restriction to one space dimension, and
the special quantities which can be analysed.
In fact,  even for the totally asymmetric simple 
exclusion process in $d=1$ there was previously no proof that
$D(t)$ diverges as $t \to \infty$.
And for general asymmetric exclusion 
processes all of these problems
were completely open. 
In this article, we present a method to study the diffusion coefficient
of general asymmetric exclusion processes.  Using this we obtain without
too much work lower bounds
$D(t) \ge (\log t)^{1/2}$ in dimension $d=2$ and
$D(t) \ge t^{1/4}$ in $d=1$. We have restricted the proof
to the case $\rho=1/2$, but a similar proof works for all 
densities away from zero or one.

\bigskip

\subsection{The model and main results}

Denote the particle configuration by $\eta = \{\eta_x\}_{x \in \Z^d}$
where $\eta_x$ is equal to $1$ if site $x$ is occupied
and is equal to $0$ otherwise. Denote by
$\eta^{x,y}$  the configuration obtained from $\eta$
by exchanging the occupation variables at $x$ and $y$~:
$$
(\eta^{x,y})_z\; =\;
\begin{cases}
\eta_z & \hbox{if $z\neq x$, $y$,}\\
\eta_x & \hbox{if $z=y$ and}\\
\eta_y & \hbox{if $z=x$.}
\end{cases}
$$
We assume that the jump law $p(\cdot)$ is local, $p(z) = 0$ for $|z|\ge L$
for some $L<\infty$, and that there is transport of mass, $\sum_z zp(z) \neq 0$.  
The  generator  of
the  asymmetric simple exclusion process is given by
\begin{equation}
({\mathcal L}f) (\eta)\; =\; \sum_{x,y\in\Z^d} p(y-x)\eta_x
(1-\eta_{y}) [f(\eta^{x,y}) -f(\eta)]\; .
\label{1g0}
\end{equation}
For each $\rho$ in $[0,1]$, denote by
$\nu_\rho$ the Bernoulli product measure on $\{0,1\}^{\mathbb Z^d}$
with density $\rho$ and by $<\cdot, \cdot>_\rho$ the inner product in
$L^2(\nu_\rho)$. The probability measures $\nu_\rho$ are invariant
for the process.

For two cylinder functions $f$, $g$ and
a density $\rho$, denote by $\la f;g\ra_\rho$ the covariance of $f$ and
$g$ with respect to $\nu_\rho$~:
$$
\la f;g\ra_\rho \; =\; \la fg\ra_\rho - \la f\ra_\rho \la g\ra_\rho.
$$
Let $P_\rho$ denote the law of the
asymmetric exclusion process starting from the equilibrium
measure $\nu_\rho$.
Expectation with respect to  $P_\rho$ is
denoted by $E_\rho$.
Let
$$
S_\rho(x, t)= E_\rho [ \{ \eta_x(t) -\eta_x(0)\} \eta_0(0)]
$$
denote the time dependent correlation functions in equilibrium
with density  $\rho$.
Denote by $\chi $  the compressibility given by
$$
\chi=\chi (\rho) = \sum_x \la \eta_x ; \eta_0 \ra_\rho.
$$
In our setting, $\chi (\rho) = \rho (1-\rho)$.

The bulk diffusion coefficient is defined by
\begin{equation}
D_{i,j}(\rho,t)  \; =\; \frac{1}{t} \frac{1}{2\chi}
\bigg\{ \sum_{x\in\Z^d}
x_i x_j S_\rho(x, t)
-\chi(v_i t)(v_j t)\bigg\},
\label{0.4}
\end{equation}
where  $v$ in $\R^d$ is the velocity defined by
\begin{equation}
vt\; =\; \frac 1{\chi}
\sum_{x\in\Z^d}  x E_\rho \Big[  \{ \eta_x(t) -\eta _x(0)\}
\eta_0(0)\Big ] \; .
\label{0.9}
\end{equation}
To simplify the notation we now specialize to  the special case
of the totally asymmetric simple exclusion process
to the right in $d=1$ and, in $d=2$, jumps only to the nearest neighbor to
the right in the $x_1$ coordinate, and jumps to both nearest neighbors in the
$x_2$ coordinate with symmetric jump rule.  More precisely, we take \begin{equation}
({\mathcal L}f) (\eta)\; =\;    \sum_{x\in\Z}  \eta_x
(1-\eta_{x+ 1 }) [f(\eta^{x, x+ 1 }) -f(\eta)], \qquad d=1; 
\label{1g}
\end{equation}
\begin{equation}
({\mathcal L}f) (\eta)\; =\;   \sum_{x\in\Z^2}  \eta_x
[1-\eta_{x+ e_1 }] [f(\eta^{x, x+ e_1 }) -f(\eta)]   +
\frac{1}{ 2} [f(\eta^{x, x + e_2 }) -f(\eta)], \qquad d=2.
\label{1g2}
\end{equation}
where we have combined the symmetric jump in the $x_2$-axis into the last term.
We emphasize that the result and method in this paper apply
to all asymmetric exclusion processes; the special choice
is made only to simplify the notation.
The velocity of the totally asymmetric simple exclusion process
is  explicitly computed as $v= 2(1-2\rho)$ in $d=1$ and 
$
v = 2 (1-2\rho) e_1$ in $d=2$.
We further assume for simplicity that the density $\rho= 1/2$
so that the velocity is zero.

Denote the instantaneous currents (i.e., the difference between
the rate at which a particle jumps from $x$ to $x+e_i$ and the
rate at which a particle jumps from $x+e_i$ to $x$) by $\tilde w_{x,x+e_i}$:
\begin{equation}
\tilde w_{x,x+ e_1 }\; =\;  \eta_x [1-\eta_{x+ e_1 }],
\qquad \tilde w_{x, x+e_2} = \frac{1}{ 2}[ \eta_{x+ e_2 } - \eta_x]
\label{0.12}
\end{equation}
We have the conservation law
$$
{\mathcal L} \eta_0 + \sum_{i=1}^d  ( \tilde w_{-e_i,0}-\tilde w_{0,e_i})=0\; .
$$
Let ${ w}_i(\eta)$ denote
the renormalized current in the $i$-th direction:
\begin{equation*}
 { w}_i(\eta)  = \tilde w_{0,e_i}  - \la \tilde w_{0,e_i}\ra_\rho  -
\frac d{d\theta}  \la \tilde w_{0,e_i}\ra_{\theta}\Big\vert_{\theta=\rho}
(\eta_0 -\rho)
\end{equation*}
Note  the subtraction of the linear term in this definition.
We have
\begin{equation}
 { w}_1(\eta)  = (\eta_0 -\rho)(\eta_{e_1 } -\rho)+\rho [\eta_{e_1} - \eta_0]
, \quad
 { w}_2(\eta) = \frac{1}{2} [ \eta_{e_2} - \eta_0 ]
\label{0.8}
\end{equation}

A function $f$ on $\{0,1\}^{\Z^d}$ will be called {\it local } if it
only depends on the variables at finitely many sites.
For local functions $f$ and $g$ we define the semi-inner product
\begin{equation}
\l  g, h \r_\rho \; =\;
 \sum_{x\in \mathbb Z^d } <  \tau_x g \, ;  \,  h >_\rho
\; =\;   \sum_{x\in \mathbb Z^d } < \tau_x h \, ;  \,  g >_\rho\; \cdot
\label{1a}
\end{equation}
Since the density $\rho$ is fixed to be $1/2$ in this article, 
we will henceforth leave out the subscript.
All but a finite number of terms in this sum vanish
because $\nu_\rho$ is a product measure and $g$, $h$ are local. 
From this inner product, we define the seminorm:
\begin{equation}
\| f \|^2 = \l  f, f \r.
\label{norm}
\end{equation}
Note that {\it gradient} terms $g= \tau_x h -h$  and 
 all degree one functions vanish in this norm. Therefore we
shall identify the currents $w$ with their degree two parts:
For the rest of the article we set
\begin{equation}
 { w}_1(\eta)  = (\eta_0 -\rho)(\eta_{e_1 } -\rho)
, \quad
 { w}_2(\eta) = 0
\label{0.8a}
\end{equation}

Fix a unit vector  $\ell \in \R^d$. A
simple calculation using Ito's formula \cite{LOY}
allows one to rewrite the diffusivity as
\begin{equation}
\ell \cdot D(t) \ell  - \frac{1}{2}
 \; = \;
\frac{1}{\chi}
\; \left\|  t^{-1/2}\int_0^t 
\ell \cdot { w} (\eta(s) )ds \right\|^2.
\label{9.3}
\end{equation}
This is a variant of the  Green-Kubo formula \cite{LOY}.
In $d=1$ of course $D(t)$ is a scalar.  In our special case in 
$d=2$, since $w_{2}=0$, $D(t)$ is a matrix with all entries zero except
$$
D_{11}(t)
 \; = \; \frac{1}{2} +
\frac{1}{\chi}
\; \left\|  t^{-1/2}\int_0^t
 { w}_1 (\eta(s) )ds \right\|^2.
$$

Recall that $\int_0^\infty e^{-\lambda t} f(t) dt \sim \lambda^{-\alpha}$
as $\lambda\to 0$ means, in some weak sense, that $f(t)\sim t^{\alpha -1}$
as $t\to \infty$.
Throughout the following $\lambda$ will always be a positive real number.
We can now state the main result.


\begin{theorem}\label{th:main}
There exists $C>0$ so that for sufficiently small
 $\lambda>0$,
\begin{eqnarray}
d=1:&& \int_0^\infty e^{-\lambda t} \; t D (t) dt \ge C \lambda^{-2-\frac{1}{4}},\\ d=2: &&
\int_0^\infty e^{-\lambda t} \; t D_{11} (t) dt \ge C \lambda^{-2}
\big | \log \lambda \, \big |^{1/2}.
\end{eqnarray}
\end{theorem}

\bigskip
The conjectured behavior for $t$ large is $D(t) \sim t^{1/3}$ in $d=1$ and 
$D_{11}(t) \sim (\log t)^{2/3}$ in $d=2$.
This theorem says that in a certain average, asymptotic  sense $D(t) \ge t^{1/4}$ in  $d=1$ and 
$D_{11}(t) \ge (\log t)^{1/2}$ in $d=2$.

From the definition, we can rewrite the diffusion coefficient as
$$
t  D_{11}(t) = \frac t 2 +  \frac 2 \chi  \int_0^t\int_0^s \;
\l e^{ u {\mathcal L}} w_1,w_1\r \; du ds
$$
in $d=2$, with an analogous formula in $d=1$ (just drop the subscripts).
Thus
\begin{eqnarray}\label{intid}
\int_0^\infty e^{-\lambda t}\;  t D_{11}(t)  dt
& = & \frac 1 {2 \lambda^2} + \frac 2 \chi \int_0^\infty dt \int_0^t\int_0^s e^{-\lambda t}\;
 \; \l e^{ u {\mathcal L}} w_1,w_1 \r \; du ds  \no \\
& = & \frac 1 {2 \lambda^2} + \frac 2 \chi \int_0^\infty  du  \Big \{
\int_u^\infty d t \;  e^{-\lambda (t-u) }
\Big ( \int_u^t d s \; \Big )\; \Big \}
 \; \l e^{-\lambda u } e^{ u {\mathcal L}} w_1 \; , \; w_1 \r  \no \\
& = & \frac 1 {2 \lambda^2} +\chi^{-1} \lambda^{-2}\l w_1
\; ,\; (\lambda-{\mathcal L})^{-1} w_1\r.
\end{eqnarray}
Therefore, Theorem \ref{th:main} follows from
the following estimate on the resolvent.

\begin{lemma}\label{le:main} There exists $C>0$ such that for sufficiently small $\lambda>0$,
\begin{eqnarray}
d=1: && \l w,(\lambda-{\mathcal L})^{-1} w\r \ge C\lambda^{-1/4};\\
d=2: && 
\l w_1,(\lambda-{\mathcal L})^{-1} w_1\r \ge C
\big | \log \lambda \, \big |^{1/2}.
\end{eqnarray}
\end{lemma}

\section{Duality and resolvent hierarchy}

Let $\calL^*$ denote the adjoint of $\calL$ with respect to the inner
product of $L^2(\nu)$ and $\calS= (\calL+\calL^*)/2$ and
$\A=(\calL-\calL^*)/2$ be the symmetric and antisymmetric parts of
$\calL$ so that 
\begin{equation} 
\calL=\calS +\A.  
\end{equation} 
Denote by $\mathcal C =
\mathcal C (\rho)$ the space of $\nu_\rho$-mean zero local functions.
For a finite subset $\Lambda$ of $\mathbb Z^d$, denote by
$\xi_\Lambda$ the $\nu_\rho$-mean zero local function defined by
\begin{equation}\label{999}
\xi_\Lambda\; =\;\prod_{x\in\Lambda} \xi_x,
\qquad \xi_x =
\frac { \eta_x-\rho} {\sqrt {\rho(1-\rho)}}\; .
\end{equation}
Note that the collection $\{\xi_\Lambda\}$ where $\Lambda$ ranges over finite
subsets of  $\mathbb Z^d$, forms an orthonormal basis of $L^2(\nu_\rho)$.
Denote by $\mathcal  M_n$
the space of local functions of degree $n$, i.e.,
the space generated by monomials of degree $n$~:
\begin{equation}\label{9999}
\mathcal  M_n\; =\; \Big\{ f\in\mathcal  C\, ;\;
f = \sum_{ |\Lambda| = n}
f_\Lambda \xi_\Lambda\, , \; f_\Lambda\in\R\Big\}\; .
\end{equation}
Note that in the definition all but a finite number of
coefficients $f_\Lambda$ vanish because $f$ is assumed to
be local. Denote by
\begin{equation}
\label{seeenn}
\mathcal  C_n=\cup_{1\le j\le n} \mathcal  M_j 
\end{equation}
 the space
of cylinder functions of degree less than or equal to $n$.
All $\nu_\rho$-mean zero local functions $f$ can be decomposed uniquely as
a finite linear combination of cylinder functions of finite
degree~: $\mathcal  C = \cup_{n\ge 1} \mathcal  M_n$.  Any $f\in L^2(\nu_\rho)$
can be written by degree,
$$
f= (f_1,f_2,f_3,\ldots)
$$
with $f_n\in \mathcal M_n$.  

For  $f\in\mathcal  M_n$ represented as
$f=\sum_{\Lambda,|\Lambda|=n} f_\Lambda \xi_\Lambda$ we have
\begin{equation}\label{funnyl}
\| f\|^2 = \sum_{n=1}^\infty \sum_{\Lambda\subset\Z^d, |\Lambda|=n} \sum_{x\in \Z^d} f_{\tau_x\Lambda} f_\Lambda.
\end{equation}
Note that this effectively reduces the degree by one.  
A simple computation shows that in this basis the symmetric part $\calS$ of
$\calL$ maps $\mathcal M_n$ into itself and is given by
\begin{equation}
\calS f(\eta) \; =\;
-\frac{1}{ 2} \sum_{j=1}^d\sum_{x\in \mathbb Z^d}  \sum_{\frac{\scriptstyle \Omega, \, |\Omega|=n-1}
{ \scriptstyle \Omega  \cap \{ x, x+e_j \} = \emptyset}}
\left[ f_{ \Omega \cup \{ x+e_j \} } -  f_{ \Omega \cup \{ x\} } \right]
\left[\xi_{\Omega \cup \{ x+e_j \}} - \xi_{\Omega \cup \{ x\}} \right].
\end{equation}
The asymmetric
part $\A$ can be decomposed into three pieces $\A= \A_0+\A_+-\A_+^*$
where $ \A_0: \mathcal  M_n\to\mathcal  M_n $ into itself, 
 $\A_+: \mathcal  M_n\to \mathcal  M_{n+1}$ and $\A_+^*: 
\mathcal  M_n\to\mathcal  M_{n-1}$ is the adjoint of $\A_+$:
\begin{align*}
&  {  \A_0} f(\eta) \; =\; \frac{1-2\rho}{2}
\sum_{x\in \Z^d}  \sum_{\frac{\scriptstyle \Omega, \, |\Omega|=n-1
}{\scriptstyle \Omega  \cap \{ x, x+ e_1 \} = \emptyset}}
\left[ f_{ \Omega \cup \{ x+ e_1 \} } -  f_{ \Omega \cup \{ x\} } \right]
\left[\xi_{\Omega \cup \{ x+ e_1 \}} + \xi_{\Omega \cup \{ x\}} \right]\; , \\
& \qquad \A_+ f (\eta) \; =\; - \sqrt{\rho(1-\rho)}
\sum_{x\in \Z^d}  \sum_{\frac{\scriptstyle \Omega, \, |\Omega|=n-1}{
 \scriptstyle \Omega  \cap \{ x, x+ e_1 \} = \emptyset}}
\left[ f_{ \Omega \cup \{ x+ e_1  \} } -  f_{ \Omega \cup \{ x\} } \right]
\xi_{\Omega \cup \{x, x+ e_1 \}} .
\end{align*}
In our special case $\rho=1/2$, we have $ \A_0=0$ and thus
$\A = \A_+ - \A_+^*$.  We can also identify $f=\sum_{\Lambda,|\Lambda|=n} f_\Lambda \xi_\Lambda$ with a symmetric function of $n$ variables,
\begin{equation}
f(x_1,\ldots,x_n) = \begin{cases} f_{\{x_1,\ldots,x_n\}}
 & {\rm if}~ x_i\neq x_j ~{\rm for}~  i\neq j,~ i,j=1,\ldots,n;\\
0 & {\rm otherwise}.\end{cases}
\end{equation}
With this notation we have 
\begin{eqnarray}\label{Afdef}
 \A_+ f(x_1, \ldots,  x_{n+1})  
 & = &-\frac  1 2 \sum_{ i = 1}^ {n+1} \sum_{j \not =  i } \delta (x_{{j}}- x_{i}-e_1 ) \prod_{k \not = j} \big (
1-\delta (x_{j}-x_{k}) \big )
\nabla_+^{ij} f  \\
\nabla_+^{ij}f(x_1, \ldots,  x_{n+1}) &=&  f(x_1, \ldots, x_{i}+ e_1,
\ldots, \widehat { x_{j}}, \ldots  x_{n+1})-
f(x_1, \ldots, x_{i},\ldots,  \widehat { x_{j}}, \ldots, x_{n+1}) \no
\end{eqnarray}
where $\delta (0)=1$ and zero otherwise, and $ \widehat { x_{j}}$ indicates the absence of $ { x_{j}}$ in
the vector. Also, 
\begin{align}\label{cals}
  \calS f(x_1, \ldots,  x_{n})
& =  \sum_{ i = 1}^ {n} \sum_{\sigma = \pm } \sum_{ j = 1}^d
 \prod_{k \not = i} \big (\, 1-\delta (x_{i} + \sigma e_j-x_{k}) \, \big )
\no \\
&   \times [f(x_1, \ldots x_{i}+ \sigma e_j, \ldots, x_{n})-
 f(x_1, \ldots, x_{i}, \ldots, x_{n})].
\end{align}
Note that $\calS$ is  the discrete Laplacian with Neumann
boundary condition on \begin{equation} \label{eone}
\E_1 = \{ \xb_n := (x_1, \ldots, x_n): x_i \not = x_j, \text{ for } i \not = j \}.
\end{equation}

The current $w\in \mathcal C_2$ 
and in our notation $w_{\{0,e_1\}}=1/4$
and $w_{\Lambda} = 0 $ for $\Lambda\neq\{0,e_1\}$, $|\Lambda|=2$.  
The resolvent equation $(\lambda-\calL) u = w$ becomes the hierarchy
\begin{align}
\A_+^* u_3+ (\lambda-\mathcal S) u_2 & = w, \\
 \A_+^* u_{k+1} + (\lambda-\mathcal S) u_k - \A_+ u_{k-1} & = 0,\qquad k\ge 3.
\end{align}
The hierarchy starts at degree $2$ instead of $1$ because the
$\l \cdot,\cdot\r$ inner product effectively reduces the degree by one.
Any term $\A_+ u_1$  is trivial in the sense of \eqref{funnyl} and
hence the degree one term plays no role and we can set $u_1=0$.
In the same way, we disregard the degree one part of the current, and
 take $w_\Lambda=0$ for all finite subsets $\Lambda\subset\Z^d$
 except $\{0,e_1\}$, and $w_{ \{0,e_1\} } = 1/4$.

Consider the truncated equation 
up to  the degree $n$,
\begin{align}\label{truncatedN}
\A_{+}^* u_{3} + (\lambda-\calS) u_2 = & w,  \no \\
\A_+^* u_{k+1} + (\lambda-{\calS}) u_k -\A_+ u_{k-1}= & 0, \quad n-1\ge 
 k \ge 3 \\
(\lambda-{\calS}) u_n -\A_+ u_{n-1}= & 0.  \no
\end{align}
We can solve  the final  equation  of \eqref{truncatedN} by
$$
 u_n = (\lambda-\calS)^{-1} \A_{+} u_{n-1}.
$$
Substituting this into the equation of degree $n-1$, we have
$$
u_{n-1} = \Big [ (\lambda-\calS)  + \A_{+}^\ast  (\lambda-\calS) ^{-1} \A_{+} \Big]^{-1} u_{n-2}.
$$
Solving iteratively we arrive at
$$
u_2 =
 \Bigg [  (\lambda-\calS) + \A_{+}^\ast \bigg \{  (\lambda-\calS)
 + \cdots +\A^*_+\Big ( (\lambda-\calS)
+ \A_{+}^\ast  (\lambda-\calS)^{-1} \A_{+} \Big )^{-1}
 \A_{+} \bigg \}^{-1}
 \A_{+} \Bigg ]^{-1} w.
$$
Let $\pi_n: \mathcal C\to \mathcal C_n$ denote the projection onto
degree $\le n$ and $\calL_n = \pi_n\calL\pi_n$. 
 The truncated equation represents the solution of
 $ (\lambda-\calL_n)u = w$
 and hence
$\l u_2,w\r = \l w  ,(\lambda-\calL_n)^{-1}w\r$. For example
\begin{align}\label{egg}
\l w, \Big [  \lambda-\calS  + \A_{+}^\ast (\lambda-\calS)^{-1}  \A_{+} \Big ]^{-1} w \r
= & \l w, (\lambda-\calL_3)^{-1} w \r. \\
\l w, \Big [  \lambda-\calS +  \A_{+}^\ast \bigg \{ \lambda-\calS+ \A_{+}^\ast
 (\lambda-\calS)^{-1}\A_{+} \bigg \}^{-1} \A_{+} \Big ]^{-1} w\r
= & \l w, (\lambda-\calL_4)^{-1} w \r.\no
\end{align}
Since $\A_{+}^\ast (\lambda-\calL)^{-1}  \A_{+}$ is nonnegative,
we have the monotonicity inequality
\begin{align}\label{mono2}
\l w, (\lambda-\calL_3)^{-1} w\r \le& \l w, (\lambda- \calL_5)^{-1} w\r\le
\cdots \le  \l w, (\lambda-\calL)^{-1} w \r\no\\ &\le \cdots \le \l w,
(\lambda-\calL_4)^{-1} w\r\le
 \l w, (\lambda-\calL_2)^{-1} w\r.
\end{align}

To check that $\l w, (\lambda-\calL)^{-1} w \r$ is in fact the limit
of these upper and lower bounds we use the variational formula.
For any matrix $M$, let $M_s$ denote the symmetric part $(M+M^*)/2$.
The identity $
\left\{  [M^{-1}]_s \right\}^{-1}
=M^\ast (M_s)^{-1} M
$ always holds,
and thus we have 
\begin{equation}\label{vf}
\l w, (\lambda-\calL)^{-1} w\r
=\sup_f\left\{ 2\l w,f\r - 
\l  (\lambda-\calL)f, (\lambda-\calS)^{-1}  (\lambda-\calL)f\r\right\}.
\end{equation}
Note that $$\l  (\lambda-\calL)f, (\lambda-\calS)^{-1}  (\lambda-\calL)f\r
=\l f, (\lambda-\calS)f\r + \l \A f, (\lambda-\calS)^{-1}  \A f\r.$$
Hence 
\begin{align}\label{vf2}
\l w, (\lambda-\calL)^{-1} w\r
=\sup_f\inf_g \left\{ 2\l w-\A^* g,f\r - \l f, (\lambda-\calS)f\r
 + \l g, (\lambda-\calS) g\r\right\}.
\end{align}
Let $a_n$ denote the supremum restricted to $f\in {\mathcal C}_n$, and
$a^n$ denote the infimum restricted to $g\in {\mathcal C}_n$ so that
$a_n\uparrow\l w, (\lambda-\calL)^{-1} w\r$ and 
 $a^n\downarrow\l w, (\lambda-\calL)^{-1} w\r$.  By straightforward
computation one checks that $a_n\le \l w, (\lambda-\calL_n)^{-1} w\r
\le a^n$, giving the desired result.

In what follows we will present a general approach to the equations (\ref{truncatedN}) which, from (\ref{mono2}) and  (\ref{egg}) give 
a nontrivial lower bound on the diffusion coefficient without too much work.
Because it gives a sequence of upper and lower bounds, the method has the
potential to give the full conjectured scaling of the diffusion coefficient.

\section{Degree 3 lower bounds}

From 
(\ref{mono2}) and  (\ref{egg}) of the previous section we have a lower
bound at degree three.  However, computations are complicated by the
hard core exclusion.  We now describe a method to remove the hard core
restriction in the computation.  We then perform the computations in 
Fourier space.  The estimates which justify the removal of the hard 
core are presented in the next section.

  Recall (\ref{eone}) the set $\E_1$
consists of configurations of $n$ particles, no two of which occupy the 
same site.  In the previous section we consider functions on $\E_1$,
and operators $\mathcal S$ and $\A_+$ acting on them.  By removal of
the hard core, we mean replacing these by functions on $\Z^{dn}$ and
operators acting on these.  
  
Suppose that $f$ is a function of particle configurations 
$f_{\{x_1, \ldots, x_n\}}$,  $(x_1, \ldots, x_n)\in \E_1$.  We can extend
$f$ to all of $ \Z^{nd}$ by defining $f (x_1, \ldots, x_n)=  f_{\{x_1, \ldots, x_n\}}$ if $(x_1, \ldots, x_n)\in \E_1$ and $f(x_1, \ldots, x_n)=0$ otherwise.
Note that
$$
E\big [ \big | \sum_{|A| = n} f_A \xi_A \big | ^2 \big ]
 = \frac 1 {n!} \sum_{x_1, \ldots, x_n \in \Z^d} |f(x_1, \ldots, x_n)|^2
$$ 
For a function $f: \Z^{nd}\to \R$, we shall use the same symbol $\langle f
\rangle $ to denote the expectation
$$
\frac 1 {n!} \sum_{x_1, \ldots, x_n \in \Z^d} f(x_1, \ldots, x_n)
$$
and write the inner product of two functions
as $\langle f,g\rangle = \langle fg \rangle$.
If $f$ and $g$ vanish on $\E_1$, this coincides with the inner
product introduced before.  We also define, as before,
$\l f, g\r = \sum_{x\in \Z^d} \langle \tau_x f, g\rangle$.
We use $\pi_n$ for the projection onto  ${\mathcal C}_n$.  

We now define $A_+F$ for symmetric functions not necessarily vanishing on
$\E_1$ by a formula analogous to \eqref{Afdef},
except that we drop the product of  delta
functions:
\begin{eqnarray}\label{AFdef}
{ A}_+ F(x_1, \ldots,  x_{n+1})   =  \frac 1 2  \sum_{ i = 1}^ {n+1}
\sum_{j \not =  i } \delta (x_{{j}}- x_{i}-e_1 )
 { \nabla}^{ij}_+ F(x_1, \ldots,  x_{n+1})
\end{eqnarray}
Notice that $\langle A_+ F \rangle = 0$ if $\langle F\rangle < \infty$
and hence the counting  measure
is invariant.  The discrete Laplacian is given by
\begin{eqnarray}\label{Deldef}
  \Delta F(x_1, \ldots,  x_{n})
 =  \sum_{ i = 1}^ {n} \sum_{\sigma = \pm } \sum_{ \alpha = 1, 2}
    [F(x_1, \ldots x_{i}+ \sigma e_\alpha, \ldots, x_{n})-
 F(x_1, \ldots, x_{i},, \ldots, x_{n})]
\end{eqnarray}
and we define
\begin{equation}
L=\Delta + A.
\end{equation}
Let $L_n = \pi_n L \pi_n$ be the restriction of $L$ to 
${\mathcal C}_n$.  

Throughout the rest of the paper we will use $C(\lambda)$ to denote a function of
$\lambda>0$ which has the property that for some $C<\infty$, and for 
sufficiently small 
$\lambda$,
\begin{equation}\label{cl}
C(\lambda)\le \begin{cases} C|\log\lambda|,& d=1; \\ C, & d=2. \end{cases}
\end{equation}
In the next section we will prove the following lemma.

\begin{lemma}\label{le9}
There exists a $C(\lambda)$ as in \eqref{cl} such that
\begin{equation}
\frac{1}{C(\lambda) n^{2}} {\l} w,(\lambda - L_n)^{-1} w {\r} \le {\l} w,(\lambda
-  {\mathcal L}_n)^{-1} w {\r} \le C(\lambda) n^2 {\l} w,(\lambda -L_n)^{-1} w {\r}.
\end{equation}
\end{lemma}
The special case $n=3$ combined with \eqref{mono2} and \eqref{egg} tells us
that
\begin{equation}
\l w, \Big [ \lambda-\Delta  + A_{+}^\ast (\lambda-\Delta)^{-1}  A_{+} \Big ]^{-1} w \r
\le C(\lambda) \l w, (\lambda -{\mathcal L})^{-1} w \r.
\end{equation}

We define the Fourier transform of $F:\Z^{nd} \to \R$ by
$$
\widehat{F} (\pb_{n}) = \sum_{\xb_n\in\Z^{nd}} e^{-i \xb_n\cdot \pb_n} F(\xb_n)
$$
for ${\pb}_n \in (\R^d/2\pi \Z^d)^n$. 
The Fourier transform of  the discrete Laplacian acting on $F$
is given by
$$
\widehat {-\Delta F} (  \pb_{n}) = - \sum_{j=1}^n \sum_{k=1}^d
\left [e^{ i  e _k  p _j}-2+  e^{- i  e _k  p _{j} } \right ] \;
\hat F( p _1, \ldots,  p _{n}) = {\omega}(\pb_n) \hat F( \pb _{n})
$$
where $ \omega( \pb_n) = \sum_{j=1}^n \omega(p_j)$.
  In $d=2$ we will denote the $d$ components of
$p$ by $(r, s)$:
$\rb_n = (r_1, \ldots, r_n)$
and $\sb_n= (s_1, \ldots, s_n)$.
 In $d=1$, $\omega(p)= -\left [e^{ i  p}-2+  e^{- i  p } \right ]$ and 
in $d=2$, $
\omega(p) = -\left [e^{ i  r}-2+  e^{- i  r } \right ]
-\left [e^{ i  s}-2+  e^{- i  s } \right ]$.  Note that in both cases
$\omega$ is real valued and nonnegative.
If $F$ is a  symmetric function of two integer variables  we have
$$
\widehat { A_+ F} ( p _1,  p _2,  p _3)
= -\frac{1}{3!} \sum_{\sigma}  \left
[e^{ i e_1 \cdot   p _{\sigma_1}}- e^{- i e_1 \cdot   p _{\sigma_3} } \right ]
\hat F( p _{\sigma_1}+ p _{\sigma_3},  p _{\sigma_2})
$$
where $\sigma$ runs over permutations of degree three.
By definition, we have
$$
\l  F,   G  \r
= \sum_z \int dp_1 d p_2 \hat F (p_1, p_2) \hat G(p_1, p_2) e^{i (p_1+p_2) z}
=  \int d p \hat F (p, -p) \hat G(p, -p)
$$
In other words, when considering the inner product $\l \cdot, \cdot \r$, we
can consider the class of $\hat F(p_1, \cdots, p_n)$ defined
only on the subspace $\sum_j p_j= 0 $  mod $2\pi \Z^d$.

\begin{lemma}\label{le4.1}
Suppose $F\in \mathcal M_2$. There exists a $C<\infty$ such that
\begin{align} d=1~: &\qquad
\l  A_+ F,  (\lambda-\Delta)^{-1} A_+   F\r\le
C \, \int_{u\in [-\pi,\pi)} 
 \omega( u) [\lambda+ \omega( u)]^{-1/2}\, |\hat F(u,-u)|^2du;\\
  d=2~: &\qquad
\l  A_+ F,  (\lambda-\Delta)^{-1} A_+   F\r\le
C \, \int_{u \in [-\pi,\pi)^2}
 \omega(e_1 \cdot   u)
  \, |\log (\lambda+  \omega(u) )| \, |\hat F(u,-u)|^2du.
\end{align}
\end{lemma}

\myproof 
Using the Schwarz inequality
to bound the cross terms by the diagonal terms, we can
bound  ${ A_+ F}$ by
\begin{align}\label{4.1}
& \l  A_+ F,  (\lambda-\Delta)^{-1} A_+   F\r = \int_{p_1+p_2+p_3=0}
 \frac {|\widehat { A_+ F} ( p _1,  p _2,  p _3)|^2}
{\lambda+ \omega(p_1)+ \omega(p_2)+\omega(p_3) } dS \no \\
& \le C     \int_{p_1+p_2+p_3=0} 
\frac {\left | e^{ i e_1 \cdot   p _1}
- e^{- i e_1 \cdot   p _3 } \right |^2
|\hat{ F}( p _1+ p _3,  p _2)|^2 }
{\lambda+ \omega(p_1)+ \omega(p_2)+\omega(p_3) }dS
\end{align}
where $dS$ is the element of surface area
 on the hyperplane $\{ p_1+p_2+p_3=0\}$.
Let $\Gamma$ denote the region in which at least one of the integration
variables, $r_i$ or $s_i$ in $d=2$ or $p_i$ in $d=1$, is bounded away from
$\pm\pi$ and $0$, let us say by $1/8$.  On $\Gamma$ the denominator,
$$
\lambda +\omega(p_1)+ \omega(p_2)+\omega(p_3)  \ge C^{-1}>0
$$
independent of $\lambda$.  Hence the integration over $\Gamma$ in
\eqref{4.1} is uniformly bounded in $\lambda$.  We are only concerned 
with terms diverging as $\lambda\downarrow 0$ and hence we can restrict
our attention to the integration over $\Gamma^C$.

We need to  divide $\Gamma^C$  according to whether $p_1, p_3$ in $d=1$ and
$r_1, r_{3},
s_1, s_{3}$ in $d=2$ are within $1/8$ of $-\pi=\pi$ or $0$. There are four
regions in $d=1$ and sixteen in $d=2$.  We have to compute each one and add
them up.  But in fact they are all analogous, and give the same result.
So for simplicity we only present the region where they are all in $[-1/8, 1/8]$.

We call $u=p_1+ p_3, v=p_1-p_3$.
The integration over corresponding region in  \eqref{4.1} is bounded
 by a constant multiple of $$
 \int_{  |u|\le 1/8}
  \omega(  u) |\hat{ F}(u, -u)|^2\left\{
 \int_{  |v| \le 1/8} \frac { dv} 
{\lambda+ |u+v|^2+ |u-v|^2+|u|^2 } \right\} du
$$
in $d=1$ and 
$$
 \int_{  |u\cdot e_1|, |u\cdot e_2| \le 1/8}
  \omega(e_1 \cdot  u) |\hat{ F}(u, -u)|^2\left\{
 \int_{  |v\cdot e_1|, |v\cdot e_2| \le 1/8} \frac { dv} 
{\lambda+ |u+v|^2+ |u-v|^2+|u|^2 } \right\} du
$$
in $d=2$.  Estimating the inside integration, in brackets,
in $d=1$,
$$
\int_{  |v| \le 1/8}  \frac { d v  }
{\lambda+ (u+v)^2+ (u-v)^2+u^2 } \le C (\lambda+ u^2)^{-1/2}.
$$
In $d=2$,
$$\int_{  |v\cdot e_1|, |v\cdot e_2| \le 1/8} \frac { dv} 
{\lambda+ |u+v|^2+ |u-v|^2+|u|^2 }
\le  C  \, |  \log (\lambda+  u^2)  | ,
$$
which completes the proof of the lemma.\myendproof

\begin{lemma}  \label{3.3}
There exists $C<\infty$ such that
\begin{align} \label{hj}
\l w, \left[ \lambda-\Delta  + A_{+}^\ast (\lambda-\Delta)^{-1}  A_{+} \right]^{-1} w \r
\ge\begin{cases} C^{-1} \lambda^{-1/4} & d=1;\\
C^{-1}  \, | \log \lambda |^{1/2} & d=2.
\end{cases}
\end{align}
\end{lemma}
\myproof
$d=1$:
The Fourier transform of the current $w$ is
$$
\hat w (p_1, p_2)= e^{-i p_2}/2
$$
From the previous lemma the left hand side 
of \eqref{hj} is bounded below by
\begin{align} &
 C^{-1} \int_{-\pi}^\pi
\frac { d \xi }  { \lambda+ \omega(\xi)
(\lambda+  \omega(\xi))^{-1/2}  }
\end{align}
We can restrict the integration to the region $\xi\in[-1/8, 1/8]$
and replace $\omega(x)$ by $x^2$ to have a further
lower bound,
$$  \int_{-1/8}^{1/8}
\frac { d \xi }  { \lambda+ \xi^2
(\lambda+  \xi^2)^{-1/2} }\ge C^{-1} \lambda^{-1/4}
$$
This proves Lemma 1.1 for $d=1$.

 $d=2$:
The Fourier transform of the current $w$ is
$$
\widehat {w} (p_1, p_2)= e^{-i r_2}/2
$$
From the previous lemma  the left hand side of
\eqref{hj} is bounded below by
\begin{align} &
C^{-1} \int_{-\pi}^\pi \int_{-\pi}^\pi
\frac { d \xi d \eta}  { \lambda+ \omega(\eta) + \omega(\xi) 
+ \omega(\xi)   |  \log (\lambda+ \omega(\eta) + \omega(\xi))  | }
\end{align}
We can restrict to the region $\xi,\eta \in [- 1 /8, 1/8]$ to have a further
lower bound. In this region, we can replace $\omega(x)$ by $x^2$ up to a
constant factor and use $| \log (\lambda+ \eta^2 +\xi^2)  |
\le | \log (\lambda+ \eta^2 ) |$ to obtain a further lower bound
\begin{equation}
 C^{-1} \int_{-1/8}^{1/8} \int_{-1/8}^{1/8} \frac { d \xi d \eta}  { \lambda+ \eta^2 \,+ \xi^2 \,
 +\, \xi^2  \, |  \log (\lambda+ \eta^2 )  | }
\end{equation}
Letting $
z = \xi (1+  \, |  \log (\lambda+ \eta^2 )  |)^{1/2}
$
the integration is bounded below by
$$
\int_{-1/8}^{1/8} \int_{-1/8}^{1/8} \frac { d z d \eta}  { \lambda+ \eta^2 \,+ z^2 \, }
(1+   \log (\lambda+ \eta^2 )  |)^{-1/2}.
$$
Changing to polar coordinates $r, \theta$ and
restricting to $ \pi /3 \le \theta \le 2 \pi /3$
we get  another lower bound,
$$
 \int_0^{1/20} \frac { rdr}  { \lambda+ r^2 \, }
 \, |  \log (\lambda+ r^2)  |^{-1/2}
 \ge  C^{-1}  \, | \log \lambda |^{1/2}  .
$$
\myendproof

Lemma \ref{le:main}  follows immediately 
from Lemma \ref{3.3}, \eqref{mono2} and Lemma \ref{le9} in $d=2$.
  However in $d=1$ this only gives the slightly
weaker result $\l w(\lambda - \calL)^{-1} \r
\ge C^{-1} \lambda^{-1/4}/ |\log\lambda|$.  
The $\log\lambda$ is purely from the removal
of the hard core, and we now show that 
the correct degree 3 lower bound in $d=1$
is of order $\lambda^{-1/4}$.

{\bf Proof of Lemma \ref{le:main} 
in $d=1$:} 
We will choose a test function $f\in \mathcal M_2$ in \eqref{mono2}
to obtain a lower bound.  Because of the summation over shifts in the
definition \eqref{1a} of $\l\cdot,\cdot\r$ there is really a reduction
in dimension from $\{ (x_1,x_2)\in\Z^2~|~ x_2>x_1\}$ to $\Z_+$.  
Hence we define for $x\in \Z_+$
$$
 \overline{f} (x) = \sum_{y\in \Z} f(y, y+1+ x).
$$
Then we have  $$2\l w,f\r=  \overline{f} (0),$$
$$
\l f,(-\calS) f\r = \frac{1}{ 2} \sum_{x\in \Z_+} [  \overline{f}  (x+1) -  \overline{f}(x) ]^2,
$$
and if $\overline{g}(x_1,x_2) =\sum_y g(y, y+1+x_1, y+2+ x_2)$,
$$
\overline{\A_+ f}(x_1,x_2) =  \delta(x_1) [\overline{ f} (x_2-1) -\overline{f}(x_2)]
+ \delta(x_2)  [\overline{f}( x_1+1) - \overline{f}(x_1)],
$$
and 
\begin{align}\label{df}
\l g, (-\calS) g\r =  
= \frac{1}{4}\sum_{x_1,x_2\in{\bf Z}_+}\sum_v
( \overline{g}((x_1,x_2) + v)- \overline{g}(x_1,x_2))^2
\end{align}
where $v$ runs over $\pm (1,0)$, $\pm (0,1)$ and $\pm (-1,1)$ with addition
intepreted as with reflecting boundary conditions at the origin.

We now make the choice  
 $$\overline{f}(x)=\lambda^{-1/4} e^{-\lambda^{3/4}x}.$$ 
 We have
$$
\lambda \sum_{x\in \Z_+} |\overline{f} (x) |^2
=\lambda^{1/2} (1-e^{-2\lambda^{3/4}})^{-1} \sim \frac{1}{2} \lambda^{-1/4}
$$
and $$ \sum_{x\in \Z_+} [\overline{f}  (x+1) - \overline{f}(x) ]^2
=
\lambda^{-1/2 }(1-e^{-\lambda^{3/4}}) \sim \lambda^{1/4},
$$
\begin{align*}
\overline{\A_+ f}(x_1,x_2)
&=-\delta(x_1){\bf 1}(x_2>0) \lambda^{-1/4}(1-e^{\lambda^{3/4}})e^{-\lambda^{3/4}x_2}
+\delta(x_2)\lambda^{-1/4} (e^{-\lambda^{3/4}} -1) e^{-\lambda^{3/4} x_1}\\
&\sim\lambda^{1/2}[\delta(x_1){\bf 1}(x_2>0) e^{-\lambda^{3/4}x_2}
-\delta(x_2) e^{-\lambda^{3/4} x_1}].
\end{align*}
and the final term $\l \A_+ f, (\lambda-\calS)^{-1} \A_+ f\r$ can be written 
$$
\int_0^\infty dt e^{-\lambda t}
\sum_{x_1,x_2,y_1,y_2\in {\bf Z}_+}
p_t((x_1,x_2), (y_1,y_2))\overline{\A_+ f}(x_1,x_2)\overline{\A_+ f}(y_1,y_2)
$$
where $p_t$ are the transition probabilities of the random walk on $\Z_+^2$
with Dirichlet form given by \eqref{df}.
Using the small $\lambda$ approximation for $\bar{\A f}$ we get
\begin{equation}\label{78}
 \lambda \int_0^\infty dt e^{-\lambda t}
\sum_{x,y\in{\bf Z}_+}e^{-\lambda^{3/4}(x+y)}\left[(1+ {\bf 1}(x,y>0)) p_t( (0,x),(0,y))
- 2{\bf 1}(x>0) p_t((0,x),(y,0))\right]
.\end{equation}
 By images we can rewrite this as
\begin{align*}
\lambda\int_0^\infty dt e^{-\lambda t} \sum_{x,y\in {\bf Z}}
e^{-\lambda^{3/4}(|x|+|y|)}\Big[ & (1+{\bf 1}(|x|, |y|>0)) p_t((0,x), (0,y))
\\ &- 2{\bf 1}(x>0) p_t((0,x),(y,0))\Big]
\end{align*}
where $p_t$ is now the transition density for a continuous time random walk on
${\bf Z}^2$ where the particle makes jumps at rate $1/2$ of $(1,0)$, $(-1,0)$,
$(0,1)$, $(0,-1)$ and, in the first and third quadrants $(1,-1)$ and $(-1,1)$,
and in the second and fourth quadrants $(1,1)$ and $(-1,-1)$.  On the axes themselves
the rules are changed a bit.  For example, on the positive $x$-axis the particle
jumps at rate $1/4$ of $(0,1)$, $(0,-1)$, $(-1,1)$ and $(-1,-1)$.  On the other axes 
these rules are just naturally rotated.  
The diffusion approximation  is
$$
2\lambda\int_0^\infty dt e^{-\lambda t} \int_{-\infty}^\infty dx\int_{-\infty}^\infty dy
e^{-\lambda^{3/4}(|x|+|y|)}[ p_t((0,x), (0,y))- p_t((0,x),(y,0))]
$$
where $p_t$ is now the transition density for a diffusion in ${\bf R}^2$
with generator $\Delta + d^2$ where $d= \partial_y-\partial_x$ in the first
and third quadrants and $d=\partial_y +\partial_x$ in the second and fourth
quadrants.  The corresponding Dirichlet form is comparable to the standard one
and therefore we can bound the transition probabilities above and
below by those of Brownian motion (see \cite{D}).
By change of variables
$$
\lambda\int_0^\infty dt e^{\lambda t}\int_{-\infty}^\infty dx\int_{-\infty}^\infty
dy e^{-\lambda^{3/4}(|x+y|)} \frac{e^{- \frac{|y-x|^2}{4t}}}{4\pi t}=C\lambda^{-1/4}.
$$
The error
$$
\lambda \int_0^\infty dt e^{\lambda t}\int_{-\infty}^\infty dx\int_{-\infty}^\infty
dy[ e^{-\lambda^{3/4}(|x|+|y|)}- e^{-\lambda^{3/4}(|x+y|)}]
 \frac{e^{- \frac{|y-x|^2}{ 4t}}}{ 4\pi t}=O(\lambda).
$$
We can choose $\alpha\phi$ for our test function instead of $\phi$.  In summary 
we have
$$
\l w,(\lambda -{\calL}^{-1}) w\r
\ge  2\alpha\lambda^{-1/4} - C\alpha^2 \lambda^{-1/4} \ge C' \lambda^{-1/4} 
$$
if  $\alpha $ is chosen sufficiently small.\myendproof

\section{Removal of hard core}

In this section we prove Lemma 3.1.  Recall the main statement
is that $\l w, (\lambda-\calL_n)^{-1} w\r$ can be bounded 
above and below in terms of $\l w, (\lambda-L_n)^{-1} w\r$
at the expense of constants depending on $n$.  Since we only use
the bound for $n=3$ in this article, the precise dependence of the
constants on $n$ is not important and in many places is probably not
optimal.

For $f$ a symmetric function on $\Z^n$,
we denote by $\rho_i , i=1, 2, 3$  the one, two and three
point functions~:
\begin{align}
 \rho_1 (x)
\; & =\; \frac 1 {(n-1)!} \sum_{x_1, \cdots,  x_{n-1}}
[f(x, x_1, \cdots, x_{n-1})]^2 \\
 \rho_2 (x,y)
\; & =\; \frac 1 {(n-2)!2!} \sum_{x_1, \cdots,  x_{n-2}}
[f(x, y, x_1, \cdots, x_{n-2})]^2 \\
 \rho_3 (x,y, z) \; & =\; \frac 1 {(n-3)!\,  3!}\sum_{x_1, \cdots,  x_{n-3}}
[f(x, y, z, x_1, \cdots, x_{n-3})]^2
\end{align}
for all $x \neq y \not = z$ in $\mathbb Z^d$. The following is
Lemma 4.8 of \cite{LY}.

\begin{lemma} \label{2.1.1}
Suppose $f$ is a symmetric function on $\Z^n$. Then we have
in all dimensions
$$
\langle \rho_3^{1/2}, (-\calS_3) \rho_3^{1/2}\rangle
\le n^2 \langle f, (-\calS) f\rangle
$$
\end{lemma}

\myproof  By symmetry, 
$$\langle \rho_3^{1/2}, (-\calS_3) \rho_3^{1/2}\rangle= 3\sum_{\sigma=\pm,~ e}
\langle (\rho^{1/2}_3( x+\sigma e,y,z)
 -\rho^{1/2}_3(x,y,z))^2\rangle
$$ 
where $e$ are summed over the standard basis of unit vectors in  $\R^n$, and 
$$ \langle f, (-\calS) f\rangle
= n\sum_{\sigma=\pm,~ e}\frac 1 {n!}\sum_{x,y,z,x_1,\ldots, x_{n-3}} [f(x+\sigma e,y,z,x_1,\ldots,x_{n-3}) - f(x,y,z,x_1,\ldots,x_{n-3})]^2.$$
 By Schwarz's inequality, $((\sum_j a_j^2)^{1/2} - (\sum_jb_j^2)^{1/2})^2
\le \sum_j (a_j-b_j)^2$, we have
\begin{align*}
&(\rho^{1/2}_3( x,y,z)
 -\rho^{1/2}_3(x',y,z))^2\\&  \le   \frac 1 {(n-3)!\,  3!}\sum_{x_1,\ldots, x_{n-3}} [f(x,y,z,x_1,\ldots,x_{n-3}) - f(x',y,z,x_1,\ldots,x_{n-3})]^2,
\end{align*}
and hence
\begin{align*}
&\langle (\rho^{1/2}_3( x+e,y,z)
 -\rho^{1/2}_3(x,y,z))^2\rangle = \frac{1}{3!} \sum_{x,y,z}
 (\rho^{1/2}_3( x,y,z)
 -\rho^{1/2}_3(x+e,y,z))^2\\
&\le\frac 1 {(n-3)!\,  (3!)^2}\sum_{x,y,z,x_1,\ldots, x_{n-3}} [f(x+e,y,z,x_1,\ldots,x_{n-3}) - f(x,y,z,x_1,\ldots,x_{n-3})]^2.
\end{align*} \myendproof

The following Lemma  is a simple extension of
Theorem 4.7 of  \cite{LY}  or Lemma 3.7 of \cite{SVY} to dimensions
$d=1, 2$. Recall the definition \eqref{eone} of $\E_1$.

\begin{lemma}\label{le:hcest}
Fix $R >0$. There exist  $C(\lambda)$  as in \eqref{cl}
such that
for  $f$ a symmetric function on $\Z^n$ vanishing on $\E_1^C$,
\begin{align*}
\sum_{i\not = j \not = k} \langle  1_{\{ |x_i-x_k|+|x_k-x_j| \le R \}}
\; f^2 \rangle \le C ( \lambda ) n^2
\langle f, (\lambda- \calS) f\rangle.
\end{align*}
\end{lemma}

\myproof 
By definition of the three point function, the left side is
a constant times
$$
\langle \,  1_{\{ |x_1-x_3|+|x_3-x_2|  \le R \}}
\, \rho_3(x_1, x_2, x_3) \, \big  \rangle
$$
By the previous Lemma, we can bound the Dirichlet form of
$g= \rho_3^{1/2}$ by that of $f$. Thus we only have to prove that
$$
\langle  1_{\{ |x_1-x_3|+|x_3-x_2| \le R \}}  g^2 \rangle \le C(\lambda) \langle
g, (\lambda-\calS) g \rangle
$$
for functions $g$ of degree three.

Recall that for configuration with three particles we have
$
\E_1 = \{ \xb_3 := (x_1, x_2, x_3): x_i \not = x_j, \text{ for } i \not = j \}
$.
We have  $g( \xb_3)=0$
whenever $\xb_3 \not \in \E_1$ and the operator $\calS$
is the discrete Laplacian on $\E_1$ with Neumann boundary conditions.
Define $ G (\xb_3)= g (\xb_3)$ if 
$\xb_3 \in \E_1$ and
$
G (\xb_3)= \Av_{\yb_3 \in \E_1, |\yb_3-\xb_3|\le 2} g(\yb_3)
$
for $\xb_3 \not \in \E_1$. 
We claim  that for $G$ so defined,
\begin{equation}\label{2.2.2}
\langle G (\lambda-\Delta) G \rangle \le C \langle g
(\lambda-\calS) g \rangle.
\end{equation}
Consider
$
\langle [G (\xb_3)- G (x_1+e_1, x_2, x_3)]^2\rangle
$
with $(x_1+e_1, x_2, x_3) \in \E_1$ and $\xb_3 \not \in\E_1$.
In this case that $G (\xb_3)$
is the average of $g(\yb_3)$ with
$|\yb_3-\xb_3|\le 2$ and $\yb_3 \in \E_1$.
We can  check that $\yb_3$ and $(x_1+e_1, x_2, x_3) \in \E_1$
can be connected via nearest-neighbor bonds in $\E_1$.
Thus we have
$
\langle [G (\xb_3)- G (x_1+e_1, x_2, x_3)]^2
\rangle \le C \langle g
(-\calS) g \rangle
$.  A
similar inequality can be checked if $(x_1+e_1, x_2, x_3) \not \in \E_1$.
and this proves \eqref{2.2.2}.

Since  $g(\xb_3)=G(\xb_3)$ on $ \E_1$
and $0$ otherwise, it is clear  that
$
\langle 1_{\{ |x_1-x_3|+|x_3-x_2| \le R \}}  g^2 \rangle
\le \langle 1_{\{ |x_1-x_3|+|x_3-x_2| \le R \}}  G^2 \rangle
$.
Thus to prove the lemma  it will suffice to prove that
$$
\langle 1_{\{ |x_1-x_3|+|x_3-x_2| \le R \}}  G^2 \rangle \le C(\lambda) \langle
G, (\lambda-\Delta) G\rangle
$$
We can drop the part of $\Delta$ in the $x_3$ direction, making the right hand
side smaller.  Hence it is enough to prove that 
$$
\langle  1_{\{ |x_1|+|x_2| \le R \}}  G^2 \rangle \le C(\lambda) \langle
G, (\lambda-\Delta) G \rangle
$$
for functions $G(x_1,x_2)$.  Call $V=  1_{\{ |x_1|+|x_2| \le R \}}$.  It is
local and bounded.  We are in $\Z^{2d}$ and we want to show that there
is a $C(\lambda)$ such that
 $
-C(\lambda) (\lambda-\Delta)\ge   V 
$ as operators,  or, equivalently
$$
V^{1/2}  (\lambda-\Delta)^{-1}V^{1/2}\le C(\lambda).
$$
Let $G_\lambda(x,y)$ be the kernel of $ (\lambda-\Delta)^{-1}$ and
 $\varphi: \Z^{2d}\to \R$.  Since $V$ is bounded by $1$,  $\sum_{x,y\in\Z^{2d} } V^{1/2}(x)
G_\lambda(x,y) V^{1/2}(y) \varphi(x)\varphi(y)$ is easily bounded
by $G_\lambda(0,0)\sum_x \varphi^2(x)$.  If $d=1$, we are in
$\Z^2$ and  $C(\lambda)=G_\lambda(0,0)\le C|\log\lambda|$.  In $d=2$,
we are in the transient case $\Z^4$ and $G_\lambda(0,0)$ is bounded
uniformly in $\lambda$ (see \cite{Spit}).
\myendproof

We now divide the complement of $\E_1$ into two sets.
We call a site $x$ an isolated double site if
$$
x_i=x_j=x, \qquad |x_k-x| \ge 5 \quad \text{for all } \; k \not = i, j.
$$  We use here the lattice distance for $x, y \in \Z^d$:
$
|x-y| = \sum_{j=1}^d |x^j-y^j|
$.
Denote by $\E_2$ the set with at most
isolated double sites and
$\E_3= [ \Z^d]^n - (\E_1\cup\E_2)$ the rest.
For a configuration
$(x_1, x_1, \cdots, x_{k}, x_{k},x_{2k+1}, x_{2k+2} \cdots, x_n)$ with
$k$  isolated double sites, we define
$F= \T  f$  by
\begin{align}
&  F (x_1, x_1, \ldots, x_{k}, x_{k},x_{2k+1}, x_{2k+2} \ldots, x_n) \no \\
 &  =  \Av_{\yb_k: |x_i-y_i| = 1 \text { for all } i}
f(x_1, y_1, \ldots, x_{k}, y_k, x_{2k+1}, \ldots, x_n)
\end{align}
If $\xb \in \E_3$, then $F(\xb)= 0$, e.g.,
$F(x, x, x+e_1, x_4, \ldots, x_n) = 0$.
We also define the restriction $\Res F$ by
$
\Res F (\xb_n) = F (\xb_n )$
 if $\xb_n \in \E_1$
and $\Res F(\xb)=0$ otherwise. Note that $\T$ and $\Res$ are not inverse to
each other although $RT$ is the identity.

\begin{lemma}\label{le2} There is a $C(\lambda)$ as in \eqref{cl} such that
for any symmetric function $f$  on  $\Z^n$ vanishing on $\E_1^C$ and 
$F= \T f$
\begin{equation}
\label{2.1}
(1/C(\lambda))  n^{-2} \langle F, (\lambda- \Delta) F\rangle \le
 \langle f, (\lambda-{\mathcal S}) f\rangle \le C \langle F, (\lambda- \Delta) F\rangle
\end{equation}
Define $\tilde F = \T \Res F$ where $\Res$ is the restriction. Then in addition,
\begin{equation}\label{eq:2.1.1}
\langle F , (\lambda- \Delta)
 F \rangle
\ge C(\lambda) n^{-2} \langle \tilde F , (\lambda- \Delta)
 \tilde F \rangle
\end{equation}
\end{lemma}

\myproof  
By definition,
$
\langle F, - \Delta F\rangle$ is given by $
\frac 1 2 \sum_{\xb, \yb : |\xb-\yb|= 1} [F(\xb)-F(\yb)]^2
$. The upper bound of  \eqref{2.1} is immediate, since $\calS$ has
Neumann boundary conditions, corresponding to dropping terms in the
Dirichlet form with either
$x$ or $y$ in $\E_1^C$.  We now prove the lower bound in  \eqref{2.1}.
We  decompose $\xb$ and $\yb$ into three sets,
$\E_1, \E_2, \E_3$, so that
$$
\langle F, - \Delta F\rangle = \sum_{\alpha=1, 2, 3;  \beta= 1, 2, 3}
\Phi_{\alpha, \beta}
$$
where
$$
\Phi_{\alpha, \beta}
= \frac 1 2 \sum_{\xb \in \E_\alpha , \yb \in \E_\beta : |\xb-\yb|= 1} [F(\xb)-F(\yb)]^2.
$$
If both  $\xb$ and $\yb$ satisfy the hard core condition,
then the contribution is
$$
\Phi_{1, 1}= \frac 1 2 \sum_{\xb, \yb : |\xb-\yb|= 1} [f(\xb)-f(\yb)]^2
\le \langle f, - \calS \, f\rangle.
$$
We can estimate terms where either $\xb$ or $\yb$ is in $\E_3$ by,
$$
|\Phi_{1, 3} + \Phi_{2, 3}+ \Phi_{3,3}|\le  C  \sum_{i\not = j \not = k} \langle
1_{\{ |x_i-x_k|+|x_k-x_j| \le R \}} \; F^2 \rangle.
$$
From the  definition of $F$, we can check that
$$
\sum_{i\not = j \not = k} \langle
1_{\{ |x_i-x_k|+|x_k-x_j| \le R \}} \; F^2 \rangle
\le C  \sum_{i\not = j \not = k} \langle
1_{\{ |x_i-x_k|+|x_k-x_j| \le R \}} \; f^2 \rangle
$$
The last term is bounded by $C(\lambda) n^2 \langle
f, (\lambda- \calS) f\rangle$ from Lemma \ref{le:hcest}. Thus we have
$$
|\Phi_{1, 3} + \Phi_{2, 3}+\Phi_{3, 3} |\le C(\lambda) n^2 \langle
f, (\lambda- \calS) f\rangle \; .
$$
We now bound  $\Phi_{1, 2}$.
In this case we have, for example, 
$
\xb = (x_1, x_1+ e_1 , x_3, \cdots, x_n)$ and $ \yb = (x_1, x_1, x_3,
\cdots, x_n)
$.
Notice that because $\xb \in \E_1$,   $\yb$ can in fact 
have at most one double site.
By assumption of isolated double sites, we have
$|x_j-x_1| \ge 5$ for all $j\ge 3$. Thus
$$
F(\yb)= \Av_{ |z-x_1|=1 } f(x_1, z, x_3, \cdots, x_n)
$$
Under the assumption $|x_j-x_1| \ge 5$ for all $j\ge 3$
we can always connect $z$ to $x_1$. By Schwarz's inequality, we then have
$$
 \sum_{x_1, x_3, \cdots, x_n} |f(x_1+ e_1, z,  x_3, \cdots,
x_n)- f(x_1, z, x_3, \cdots, x_n)|^2
\le  C \langle \, f, (-\calS) f\, \rangle \,
$$
Hence
$$
|\Phi_{1, 2} | \le C \langle \, f, (-\calS) f\, \rangle \,.
$$
Finally we bound  $\Phi_{2, 2}$. The typical case looks like
$
\xb = (x_1, x_1, x_3, x_3+ e_1 , x_5,\cdots, x_n)$ and $
\yb = (x_1, x_1, x_3, x_3, x_5,\cdots,  x_n)
$.
Then
\begin{align*}
& F(x_1, x_1, x_3, x_3+ e_1 , x_5,\cdots,  x_n)
- F(x_1, x_1, x_3, x_3, x_5,\cdots,  x_n)\\
& = \Av_{ |z-x_1|=1 }
 \Av_{ |w-x_3|=1 } [f(x_1, z, x_3, x_3+ e_1 , x_5, \cdots, x_n)
-f(x_1, z, x_3, w, x_5, \cdots, x_n) ]
\end{align*}
Using Jensen's inequality and the same arguments as in the estimate of $\Phi_{1, 2}$ above we obtain $|\Phi_{2,2}| \le C \langle  f, (-\calS) f \rangle$.
Putting all these estimates together, we have the lower bound  of \eqref{2.1}.

To prove \eqref{eq:2.1.1}, call
  $f = \Res F$ so that $\tilde F = \T f$. From \eqref{2.1}
we have
$
\langle f ,
(\lambda- \calS)
 f \rangle \ge C(\lambda) n^{-2} \langle \tilde F , (\lambda- \Delta) \tilde F \rangle 
$.
Since $\calS$ is an operator with Neumann boundary condition, for
any $F$ with $\Res F = f$ we have
$
\langle F , (\lambda- \Delta)
 F \rangle
\ge C \langle f , (\lambda- \calS)
 f \rangle
$ and this proves \eqref{eq:2.1.1}.
\myendproof

\begin{lemma}\label{le3}
Suppose that $f, g\in \mathcal C_n$. Let
$F, G = \T f, \T g$. There is a  $C(\lambda)$ as in \eqref{cl} such that
\begin{equation}\label{2.2}
\big | \langle g, \A f \rangle- \langle G, A F \rangle \big | \le
C(\lambda) \; n^{1/2}\; \langle G,(\lambda- \Delta) G \rangle^{1/2}\;
 \langle F , (- \Delta)
 F \rangle^{1/2} .
\end{equation}
Recall $\tilde F = \T \Res F$. We
have
\begin{equation}\label{2.8}
\Big | \langle  A \tilde F, G \rangle - \langle  A  F, G\rangle
\Big | \le  C(\lambda) n^{1/2}
\langle G,(\lambda- \Delta) G \rangle^{1/2}
 \langle F , (- \Delta)
 F \rangle^{1/2}
\end{equation}
where $C$ is independent of $\lambda$.
\end{lemma}

\myproof We first prove \eqref{2.2}. Note that it suffices to prove it
 with $\A$ replaced by $\A_+$, since it then follows for $\A_+^*$ and
thus for
 $\A= \A_+-\A_+^*$.  Suppose first of all that
$f\in \mathcal M_n$.
Recall the definition of $\A_+$ in \eqref{Afdef}.
It suffices to consider one term in the summation, say
$\A_+^{ij} f= \delta(x_2-x_1-e_1)\prod_{k\neq 2} (1-\delta(x_2-x_k)) \nabla_+^{1,2} f$ and $A_+^{ij} f= \delta(x_2-x_1-e_1) \nabla_+^{1,2} f$.  We divide $\Z^{(n+1)d}$ into  $
{\cal N}=\{ |x_j -x_1| >  5 \text { for all } j\ge 3 \} .
$ and its complement.
If $\xb_{n+1}\in  {\cal N}$, then
$ ({ A}_+^{12} F - { \A}_+^{12} f) (\xb_{n+1}) = 0$ so we can estimate
\begin{align*}
 \left|  \langle  g, 1_{{\cal N}^C}\A^{12}_+ f  
\rangle  \right|^2
 \le \big \|\; f(x_1+ e_1 , x_3 \cdots,  x_{n+1})-
f(x_1, x_3\cdots,  x_{n+1}) \; \big \|^2 \,
 \| g 1_{{\cal N}^C} \|^2
\end{align*}
Clearly,
$
1_{{\cal N}^C} \le \sum_{j\ge 3 }  \{ |x_j -x_1|\le 5  \}
$.
By Lemma \ref{le:hcest}, 
\begin{equation} 
 \| g 1_{{\cal N}^C} \|^2 \le
C (\lambda) n^2   \langle g,(\lambda- \calS) g \rangle.
\end{equation}
Replacing $1, 2$ by $i, j$ and summing over all $i, j$, by the permutation 
symmetry of $f$  we have
\begin{equation}\label{2.21}
\big | \langle g,1_{ {\cal N}^C} \A_+ f \rangle |^2 \le
C(\lambda)   n^{2}  \langle g,(\lambda- \calS) g \rangle
  n^{-1}\langle f,(- \calS) f \rangle
\end{equation}
A similar bound holds for $F$ on ${\cal N}^c$. Combining these estimates
and using Lemma \ref{le2}, we obtain \eqref{2.2} for
$\A_+$.

If $f\in \mathcal C_n$, write $f=(f_1,\ldots,f_n, 0,0,\ldots)$.  From 
\eqref{2.21} we have
$$
\big | \langle g_{(k+1)},  \A_+ f_{(k)} \rangle | \le
C(\lambda)   n^{1/2}  \e \langle g_{(k+1)},(\lambda- \calS) g_{(k+1)} \rangle
 +  C   n^{1/2} \e^{-1}  \langle f_{(k)},(- \calS) f_{(k)} \rangle
$$
Summing over $k$ and optimizing $\e$, we get \eqref{2.2}.  Repeating the
proof, i.e., noting that
$A_+^{12}F-A_+^{12}\tilde F=0$ on $\mathcal N$ and using Lemma \ref{le:hcest}, to bound the term on $\mathcal N^C$,
gives \eqref{2.8}.
\myendproof

\begin{lemma}\label{le4}  There is a $C(\lambda)$ as in \eqref{cl}
 such that for all $f\in \mathcal C_n$, if
$F = \T f$
\begin{align} \label{2.3}
& C^{-1} \langle  A F,  (\lambda -\Delta)^{-1} A F\rangle - C^2(\lambda)n 
\langle   F,  (\lambda- \Delta) F\rangle  \no \\
& \quad\le \langle \A f,  (\lambda- \calS)^{-1} \A f\rangle  \no \\
&  \quad\quad\le  C(\lambda) n^2
\langle  A F,  (\lambda -\Delta)^{-1} A F\rangle + C^3(\lambda) n^3
\langle   F,  (\lambda- \Delta) F\rangle\; .
\end{align}
Recall $\tilde F = \T \Res F$. Then we have
\begin{equation} \label{2.3.1}
\langle   A F,  (\lambda- \Delta)^{-1} A F \rangle \ge C^{-1} \langle  A
\tilde F  (\lambda- \Delta)^{-1} A \tilde F   \rangle - C^2(\lambda) n
\langle F , (\lambda- \Delta)
 F \rangle
\end{equation}
\end{lemma}

\myproof
Recall the variational formula
$$
\langle \A f,  (\lambda- \calS)^{-1} \A f\rangle = \sup_g  \left\{
2 \langle   g, \A f \rangle - \langle   g,  (\lambda- {\calS}) g\rangle
\right\}
$$
By Lemmas \ref{le2} and \ref{le3},  we can bound the right side from above by
\begin{align*}
& 2 \langle   g, \A f \rangle - \langle   g,  (\lambda- {\calS}) g\rangle \\
&\le
  2 \langle   G, A F \rangle - (1/C(\lambda)) n^{-2}  \langle   G,  (\lambda-
\Delta) G\rangle + C(\lambda)n^{1/2} \langle G,(\lambda- \Delta) G \rangle^{1/2}
 \langle F , (\lambda-\Delta)
 F \rangle^{1/2} \\
& \le C(\lambda) n^2  \langle  A F,  (\lambda-\Delta)^{-1} A F\rangle + C^3(\lambda)  n^3
\langle   F,  (\lambda-\Delta) F\rangle
\end{align*}
which give the upper bound in \eqref{2.3}. 
Alternatively we can bound the right side below by
\begin{align*}
& 2 \langle   g, \A f \rangle - \langle   g,  (\lambda- {\calS}) g\rangle \\
&\ge   2 \langle   G, A F \rangle - C(\lambda) n^{1/2} \langle G,(\lambda- \Delta) G
\rangle^{1/2}
 \langle F , (\lambda-\Delta)
 F \rangle^{1/2}-
C \langle   G,  (\lambda- \Delta) G\rangle \\
 &\ge   2 \langle   G, A F \rangle -  C  \langle   G,  (\lambda-
\Delta) G\rangle
 - C^2(\lambda)n  \langle   F,  (\lambda-\Delta) F\rangle
\end{align*}
Optimizing over $G$, 
\begin{align*}
 2 \langle   g, \A f \rangle - \langle   g,  (\lambda- {\calS}) g\rangle 
\ge   C \langle  A F,  (\lambda-\Delta)^{-1} A F\rangle - C^2(\lambda)n  \langle   F,  (\lambda-\Delta) F\rangle
\end{align*}
This proves the lower bound in \eqref{2.3}.  By definition,
\begin{equation}\label{ASA}
\langle  A F,  (\lambda- \Delta)^{-1} A F\rangle =
\sup_G \left\{ 2\langle AF, G\rangle - \langle
G (\lambda - \Delta) G \rangle \right\}
\end{equation}
By Lemma \ref{le3},
the right hand side of \eqref{ASA} is bounded by
\begin{align*}
& 2\langle AF, G\rangle - \langle
G (\lambda - \Delta) G \rangle  \\
&\ge 
2\langle  A \tilde F, G \rangle - C(\lambda) n^{1/2} \langle G,(\lambda- \Delta)
G \rangle^{1/2}
 \langle F , (\lambda- \Delta)
 F \rangle^{1/2} -
\langle    G,  (\lambda- \Delta) G \rangle
\\
&\ge 
2 \langle  A \tilde F, G \rangle - C \langle G,(\lambda-
\Delta) G \rangle - C^2(\lambda) n \langle F , (\lambda- \Delta)
 F \rangle.
\end{align*}
Taking the sup over $G$, we prove \eqref{2.3.1}.
 \myendproof

Our interest is in the inner product $\l, \r$. It is easy to check
that all previous Lemmas hold for the inner
product $\l, \r$ as well. Since extensions from $\langle, \rangle$
to $\l, \r$ were carried out in detail in \cite{LY},
we shall only outline the basic procedures to
prove these Lemmas for the inner product $\l, \r$.

We now prove the analogue of
Lemma \ref{le2} for the inner product $\l, \r$.
For any local functions we  rewrite the inner product as
\begin{equation}
{\l}  g, h {\r} \; =\;
\lim_{k \to \infty} (2k+1)^{-2}  \langle \sum_{ |x| \le k}
\tau_x g \, ;  \,
 \sum_{ |x| \le k}   \tau_x h \rangle
\end{equation}
Similarly, we have
\begin{equation}
{\l}  g, (- \calS) h {\r} \; =\;
\lim_{k \to \infty} (2k+1)^{-2}  \langle \sum_{ |x| \le k}
\tau_x g \, ;  \,
(- \calS) \sum_{ |x| \le k}   \tau_x h \rangle
\end{equation}

Recall the definition of $F= \T  f$ is linear in $f$.
Furthermore, the intersection properties of
$(x_1, \cdots, x_n)$ are independent
of the translation, i.e.\ $\xb \in \E_i $ iff $\tau_z\xb \in \E_i$.
Thus we have
$$
\T \sum_{ |x| \le k}   \tau_x h
= \sum_{ |x| \le k}   \tau_x \T h
$$
For each fixed $k$ everything is still local so by Lemma \ref{le2}, 
\begin{align*}
 (2k+1)^{-2} \langle \sum_{ |x| \le k}
\tau_x  F, (\lambda- \Delta) \sum_{ |x| \le k}
\tau_x  F\rangle  \le
C(\lambda)  n^2 (2k+1)^{-2} \langle \sum_{ |x| \le k}
\tau_x  f, (\lambda-
{\mathcal S}) \sum_{ |x| \le k}
\tau_x  f\rangle
\end{align*}
Lettting $k \to \infty$ limit and using that the limits exist on both side
we obtain Lemma \ref{le2} for $\l\cdot,\cdot\r$.
To prove lemmas \ref{le3} and \ref{le4} for $\l\cdot,\cdot\r$, 
we only have to use that
$A$ commutes with translations.

{\bf Proof of Lemma 3.1} {\it Upper Bound}: 
We start with the variational formula (see \eqref{vf}),
\begin{equation}\label{wlow}
 \l  w,(\lambda -{\mathcal L_n})^{-1} w{\r}  =
\sup_{f\in \mathcal C_n } \left\{ 2{\l} w, f{\r} -
{\l} (\lambda-{\mathcal L_n}) f, (\lambda- {\mathcal S})^{-1}
(\lambda - {\mathcal L_n})f{\r}\right\}
\end{equation}
where $\mathcal C_n$ denotes the functions of degree less than or equal
to $n$.
By definition,
$$
{\l} (\lambda-{\mathcal L_n}) f, (\lambda-{\mathcal S})^{-1}  (\lambda -
{\mathcal L}_n)f {\r} =
{\l} f, (\lambda-\calS)  f {\r} +{\l} \A f,
(\lambda-\calS)^{-1}   \A f {\r}
$$
Here  we have
set ${\mathcal A}_{+}=0 $ on $\mathcal M_n$
since we are considering only ${\mathcal L}_n$.
By  Lemmas \ref{le2} and \ref{le4},
\begin{align*}
{\l} f, (\lambda-\calS)  f {\r} +{\l} \A f,
(\lambda-\calS)^{-1}   \A f {\r}
\ge & ( n^{-2}/C(\lambda))  {\l} F, (\lambda- \Delta) F{\r}
\\ & +n^{-2} C(\lambda)  {\l} A
F,  (\lambda-\Delta)^{-1} A F{\r}.
\end{align*}
Since $ \, {\l} f, w{\r}$  and 
${\l} F, w {\r} \, $ 
vanish except for functions of degree 2, 
it is easy to check that
$ \, {\l} f, w{\r}= {\l} F, w {\r} \, $.  Thus we
have
\begin{align*}
& {\l} w,(\lambda -{\mathcal L}_n)^{-1} w {\r} \\
& \le   \sup_{F = \T f} \left\{
 2{\l} F, w {\r}
-(1/ C(\lambda)) n^{-2}  {\l}  A F,  (\lambda-\Delta)^{-1} A F{\r}
- C(\lambda) n^{-2} {\l}   F,  (\lambda- \Delta) F{\r} \right\} \\
& \le C(\lambda) n^2  {\l} w,(\lambda -L_n)^{-1} w {\r}.
\end{align*}

\noindent {\it Lower Bound}: We have
\begin{align*}
{\l} w,(\lambda -L_n)^{-1} w {\r} = & \sup_{F}  \left\{
2 {\l} F, w {\r}
-   {\l}  A F,  (\lambda-\Delta)^{-1} A F{\r} -   {\l}   F,
(\lambda- \Delta) F{\r} \right\}
\end{align*}
Since ${\l} \tilde F, w {\r} = {\l}  F, w {\r} $,
from Lemma \ref{le4} we have
$$
{\l} w,(\lambda -L_n)^{-1} w {\r} \le \sup_F \left\{ 2{\l} \tilde F, w {\r} -
C^{-1} {\l} F , (\lambda- \Delta)
 F {\r}- \frac{n^{-1}}{ C^2(\lambda)}  {\l} A \tilde F, \;   (\lambda- \Delta)^{-1} A \tilde F   {\r}  
\right\}.
$$
By Lemma \ref{le2}, 
$
{\l} F , (\lambda- \Delta)
 F {\r}
\ge C(\lambda)  n^{-2} {\l} \tilde F , (\lambda- \Delta)
 \tilde F {\r}
$.
Recall $\Res F = f, \tilde F = \T f$. Thus the previous line is bounded above
by
\begin{align*}
 \sup_{f = \Res F }  \left\{  2 {\l}  \T f, w  {\r}
-  \frac{1}{ C(\lambda)} n^{-2}
   {\l} \T f , (\lambda- \Delta)
\T f {\r}   -  \frac{1}{ C^2(\lambda)}n^{-1}  {\l}  A \T f , \;  (\lambda- \Delta)^{-1} A
\T f {\r}
 \right\}
\end{align*}
Clearly, ${\l}  \T f, w  {\r} = {\l}  f, w  {\r} $.
By Lemma \ref{le2}, we have
$
{\l} \T f , (\lambda- \Delta) \T f {\r} \ge C {\l}  f ,
(\lambda- {\mathcal S}) f {\r}.
$
By Lemma \ref{le4}, we have
\begin{align*}
{\l} w,(\lambda -L_n)^{-1} w {\r}
& \le
2 {\l}  f, w  {\r} -  C^{-1} n^{-2}     {\l} f ,
(\lambda- \calS) f {\r}  - (1/C(\lambda)) n  {\l}  \A  f , \;
(\lambda- \calS)^{-1} \A
f {\r}  \\
& \le C(\lambda) n^{2}  {\l}  w,(\lambda - {\mathcal L}_N)^{-1} w {\r}
\end{align*}
This proves the Lemma.
\myendproof

\noindent{\bf Acknowledgement}:
H.T.Yau would like to thank  P. Deift, J. Baik
and H. Spohn for explaining their results.
In particular, Spohn has pointed out the relation
\eqref{0.1} so that the connection between the current
across the zero and the diffusion coefficient becomes
transparent. H.T.Yau would also like to thank A. Sznitman for
his hospitality and invitation to lecture on this
subject at ETH.

\end{document}